\documentclass[12pt]{article}
\newtheorem{teo}{Theorem}[section]
\newtheorem{cor}[teo]{Corollary} 
 \newtheorem{rem}[teo]{Remark}

  \def\halmos{{\rule{0.2cm}{0.2cm}}}
\def\DZ{{Z\kern-.40em Z}}
\def\DN{{I\kern-.13em  N}} \def\DH{{I\kern-.13em  H}}
\def\DR{{I\kern-.13em R}}
\def\DP{{I\kern-.13em P}}               \def\DK{{I\kern-.13em K}}
\def\DC{{\kern.24em \vrule width.05em height1.5ex depth-.05ex
\kern-.30em C}}
\def\DO{{\kern.24em \vrule width.05em height1.5ex depth-.05ex
\kern-.30em O}} \def\DQ{{\kern.24em \vrule width.05em height1.5ex
depth-.05ex \kern-.30em Q}}
     \def\dz{\mbox{\scriptsize{Z\kern-.30em Z}}}
 
 \def\dR{\mbox{\sf\DR}}

\def\dC{\mbox{\sf\DC}}

\def\Im{\mathop{\mathrm{Im}}}

\begin{document}
 \title{An explicit formula for the matrix logarithm}

\author{Jo\~ao R. Cardoso\thanks{Work
supported in part  by ISR and a PRODEP grant -- Program n.
4/5.3/PRODEP/2000.}\\ Instituto Superior de Engenharia de
Coimbra\\ Quinta da Nora\\
 3040-228 Coimbra -- Portugal\\ \texttt{jocar@isec.pt}}

 \maketitle
\begin{abstract}
We present an explicit polynomial formula for evaluating the
principal logarithm of all matrices lying on the line segment
$\{I(1-t)+At:t\in [0,1]\}$ joining the identity matrix $I$ (at
$t=0$) to any real matrix $A$ (at $t=1$) having no eigenvalues on
the closed negative real axis. This extends to the matrix
logarithm the well known Putzer's method for evaluating the matrix
exponential.

\vspace*{0.2 cm} \noindent \small{{\bf Key-words:} Matrix
logarithm, polynomial method}

\vspace*{0.2 cm} \noindent\small{{\bf AMS subject
classifications:} 15A18, 34A30}
\end{abstract}

\section{Introduction}
Given a nonsingular matrix $A\in\dR^{n\times n}$, any solution of
the matrix equation $e^X=A$, where $e^X$ denotes the exponential
of the matrix $X$, is called {\em logarithm} of $A$. In general, a
nonsingular real matrix may have an infinite number of real and
complex logarithms. If $A$ has no eigenvalues on the closed
negative real axis then $A$ has a unique real logarithm with
eigenvalues in the open strip $\{ z\in \dC :-\pi <\Im z<\pi\}$ of
the complex plane (see, for instance, \cite{Horn94}). This unique
logarithm may be written as a polynomial in $A$ and is called the
{\em principal} logarithm of $A$. It will be denoted by $\log A$.

The problem of computing the principal matrix logarithm has
received some attention in recent years (see, for instance,
\cite{Cardoso01}, \cite{Kenney01}, \cite{Dieci96}, \cite{Dieci00}
and \cite{Kenney98}). In part, this interest has been motivated by
the applications of the matrix logarithm in areas such as Systems
Theory and Control Theory. The above cited papers list some
applications.

As far as we know, most of the methods proposed for computing the
principal logarithm are approximation methods . Unlike the matrix
exponential case, for which several closed forms based on
polynomial representations have been studied (see, for instance,
\cite{Kirchner67}, \cite{Leonard96}, \cite{Putzer66} and
\cite{Leite99}), little attention has been paid to closed forms
for the matrix logarithm.

In this paper, we find for the matrix logarithm the analogue of
the well known Putzer's method \cite{Putzer66} for evaluating the
matrix exponential. Assuming that for $t\in\dR$ the spectrum of
$I-At$ does not intersect $\dR^-_0$, we consider the curve
$t\mapsto\log(I-At)$ in $\dR^{n\times n}$. Using the coefficients
of a polynomial $p(\lambda)$ of degree $k$ such that $p(A)=0$,
every matrix in that curve will be written as a linear combination
of the matrices $I,A,\cdots,A^{k-1}$, in the following way:
$$\log(I-At)=f_1(t)I+f_2(t)A+\cdots +f_k(t)A^{k-1},$$ where the
coefficients $f_1,\cdots,f_k$ are integrals of certain rational
functions.

We find this simple method suitable for teaching purposes because
the topics required for understanding it (basically, eigenvalues
of matrices and integration of rational functions) are usually
taught in the first years of undergraduate courses. We recall
that, in contrast, other methods proposed for evaluating the
matrix logarithm require advanced theory, such as Schur
decompositions, matrix square roots and matrix Pad\'{e}
approximants.

\section{A polynomial formula for the matrix logarithm}

Given $A\in\dR^{n\times n}$, let
$p(\lambda)=\lambda^k+c_1\lambda^{k-1}+\cdots +c_{k-1}\lambda+c_k$
be a polynomial with real coefficients such that $p(A)=0$ and let
$$C=\left[
\begin{array}{ccccc}
0&0&\cdots&0&-c_k\\ &&&&-c_{k-1}\\&I_{k-1}&&&\vdots\\&&&&-c_1
\end{array}
\right],$$ where $I_m$ denotes the $m\times m$ identity matrix, be
the companion matrix of $p(\lambda)$. Examples of polynomials
$p(\lambda)$ such that $p(A)=0$  are the characteristic polynomial
of $A$ ($k=n$) and the minimum polynomial of $A$ ($k\leq n$).

Before stating our main result, let us define the following subset
of $\dR$:

$${\cal D}=\left\{
t\in\dR:\sigma\left(I-At\right)\cap\dR_0^-=\phi\right\},$$ where
$\sigma(X)$ stands for the spectrum of $X$ and $A$ is a given
$n\times n$ matrix. For each $t\in\dR$, the eigenvalues of $I-At$
are of the form $1-\lambda t$, with $\lambda\in\sigma(A)$. Since
non real eigenvalues of $A$ always give rise to non real
eigenvalues of $I-At$, it is enough to consider real eigenvalues
of $A$ to obtain a more clear description of the set ${\cal D}$.
Thus, we may write $${\cal D}=\{t\in\dR :1-\lambda t>0,\quad
\forall \lambda\in\sigma(A)\cap\dR\}.$$

Let $\lambda_M=\max (\sigma(A)\cap \dR)$ and $\lambda_m=\min
(\sigma(A)\cap \dR)$. Assuming that $A$ has both positive and
negative real eigenvalues, we have ${\cal D}=]\frac{1}{\lambda_m}
,\frac{1}{\lambda_M}[$. If $A$ does not have negative eigenvalues
then ${\cal D}=]-\infty ,\frac{1}{\lambda_M}[$ and if $A$ does not
have positive eigenvalues then ${\cal D}=] \frac{1}{\lambda_{m}},
+\infty [$. In any case, ${\cal D}$ is an open interval.

\begin{teo}
Suppose that the above notation holds and that the vector function
$[f_1(t),\cdots,f_k(t)]^T$ is the solution in ${\cal D}$ of the
initial value problem \begin{equation} \label{IVP1}
\left(I-Ct\right)\dot{x}(t)=-e_2 ,\quad x(0)=0,
\end{equation}
where $e_2=[0\; 1\; 0 \cdots 0]^T$. Then
\begin{equation} \label{formula}
\log\left(I-At\right)=f_1(t)I+f_2(t)A+\cdots+f_k(t)A^{k-1},
\end{equation}
for all $t\in {\cal D}$.
\end{teo}

\noindent {\bf Proof.} The function $X(t)=\log\left(I-At\right)$
is differentiable for all $t\in {\cal D}$ and
$\dot{X}(t)=-A(I-At)^{-1}$ (see (6.6.14) and (6.6.19) in
\cite{Horn94}). Besides, $X(t)$ is the unique solution in ${\cal
D}$ of the initial value problem
\begin{equation} \label{IVP2}
\left(I-At\right)\dot{Y}(t)=-A,\quad Y(0)=0,
\end{equation}
where $Y(t)\in\dR^{n\times n}$.

Let $[f_1(t),\cdots,f_k(t)]^T$ be the solution  of (\ref{IVP1}) in
${\cal D}$ and define
$P(t):=f_1(t)I+f_2(t)A+\cdots+f_k(t)A^{k-1}$. In the following we
show that $P(t)$ is also a solution of (\ref{IVP2}). Clearly
$P(0)=0$ since $f_i(0)=0,\; \forall i=1,\cdots,k$.

Since the vector function $[f_1 \cdots f_k]^T$ satisfies
(\ref{IVP1}), a little calculation lead us to the system

\begin{equation}\label{sistema1}
\left\{\begin{array}{lcl}
 \dot{f}_1+c_kt \dot{f}_k&=&0  \\-t
 \dot{f}_1+\dot{f}_2+c_{k-1}t\dot{f}_k&=&-1\\
-t\dot{f}_2+\dot{f}_3+c_{k-2}t\dot{f}_k&=&0\\
 \cdots&&\\-t
 \dot{f}_{k-1}+(1+c_1t)\dot{f}_k&=&0
\end{array}\right..
\end{equation} Using the equations of (\ref{sistema1}) and the
identity $A^k=-c_1A^{k-1}-\cdots -c_{k-1}A-c_kI$, which follows
from the Cayley-Hamilton theorem, we may write

$$\begin{array}{rcl} (I-At)\dot{P}(t)&=&(I-At)( \dot{f}_1
I+\dot{f}_2A+\cdots+ \dot{f}_kA^{k-1})\\&=&
\dot{f}_1I+(\dot{f}_2-t\dot{f}_1)A+\cdots+(\dot{f}_k-t
\dot{f}_{k-1})A^{k-1}-t \dot{f}_kA^k\\&=&( \dot{f}_1+c_kt
\dot{f}_k)I+(-t
 \dot{f}_1+\dot{f}_2+c_{k-1}t\dot{f}_k)A+\cdots\\&&\cdots+(-t
 \dot{f}_{k-1}+(1+c_1t)\dot{f}_k)A^{k-1}\\&=&-A
\end{array}.$$

\noindent Since (\ref{IVP2}) has a unique solution, it follows
that $P(t)=\log\left(I-At\right)$. $\halmos$

\vspace*{0.2cm} Since the coefficients functions in
(\ref{formula}) are solutions of (\ref{sistema1}), we can obtain
formulae for $\dot{f}_i,\; i=1,\cdots,k$, by solving the first
equation for $\dot{f}_1$ and substituting it into the second
equation, solving the second equation for $\dot{f}_2$ and
substituting it into the third equation and proceeding similarly
until the last equation. The result is $$\begin{array}{rcl}
\dot{f}_1&=&-c_k\dot{f}_kt\\ \dot{f}_i&=&-t^{i-2} -\dot{f}_k
\sum_{j=1}^{i}c_{k-i+j}t^j, \quad i=2,\cdots,k-1,\\
\dot{f}_k&=&\frac{-t^{k-2}}{1+c_1t+\cdots +c_kt^k},
\end{array}$$
or, equivalently, \begin{center} $\begin{array}{rcl}
\dot{f}_1&=&\frac{c_kt^{k-1}}{1+c_1t+\cdots +c_kt^k}\\ &&\\
\dot{f}_i&=&\frac{-t^{i-2}-c_1t^{i-1}-
\cdots-c_{k-i}t^{k-i}}{1+c_1t+\cdots +c_kt^k}, \quad
i=2,\cdots,k-1\\&&\\ \dot{f}_k&=&\frac{-t^{k-2}}{1+c_1t+\cdots
+c_kt^k}.
\end{array}$\end{center}

We note that the constants arising in the integration process to
find $f_i,\ i=1,\cdots,k$, can be evaluated according to the
identities $f_i(0)=0,\; i=1,\cdots,k$.

\vspace*{0.3cm} We now summarize the previous discussion in the
next corollary.

\begin{cor}
Given $A\in\dR^{n\times n}$, let
$p(\lambda)=\lambda^k+c_1\lambda^{k-1}+\cdots +c_{k-1}\lambda+c_k$
be a polynomial with real coefficients such that $p(A)=0$. If
${\cal D}=\{ t\in\dR:\sigma\left(I-At\right)\cap\dR_0^-=\phi\}$,
then, for all $t\in {\cal D}$,
$$\log\left(I-At\right)=f_1(t)I+f_2(t)A+\cdots+f_k(t)A^{k-1},$$
where $f_1,\cdots,f_k$ are differentiable functions in ${\cal D}$
given by
\begin{equation}\label{formula2}
\begin{array}{rcl}
f_1(t)&=&\int_0^t\frac{c_ks^{k-1}}{1+c_1s+\cdots +c_ks^k}\ ds\\
&&\\ f_i(t)&=&\int_0^t\frac{-s^{i-2}-c_1s^{i-1}-
\cdots-c_{k-i}s^{k-i}}{1+c_1s+\cdots +c_ks^k}\ ds, \quad
i=2,\cdots,k-1\\&&\\
f_k(t)&=&\int_0^t\frac{-s^{k-2}}{1+c_1s+\cdots +c_ks^k}\ ds.
\end{array}
\end{equation}

\end{cor}

\begin{rem} {\rm There exists a relationship between the
polynomial $p(\lambda)$ and the polynomial
$q(\lambda)=1+c_1\lambda+\cdots+c_k\lambda^k$ in the denominator
of the functions under integral symbols in (\ref{formula2}):
$q(\lambda)=\lambda^k p(1/\lambda)$.}
\end{rem}

\begin{rem} {\rm The indefinite integrals in (\ref{formula2}) may be
obtained explicitly because we are dealing with rational
functions. We note that many calculus textbooks provide methods
for evaluating integrals of these kind of functions. Also,
symbolic software packages like Mathematica, Maple or Derive are
able to compute them. }
\end{rem}

\vspace*{0.3cm} \begin{rem} {\rm For $A$ such that
$\sigma(A)\cap\dR^-_0=\phi$, equation (\ref{formula}) allows us to
find an explicit formula for evaluating the logarithm of all
matrices on the line segment joining $I$ (at $t=0$) to $A$ (at
$t=1$): $$\{I(1-t)+At:t\in[0,1]\}.$$ Indeed, $$\begin{array}{rcl}
\log\left(I(1-t)+At\right)&=&\log\left(I-(I-A)t\right)\\
&=&f_1(t)I+f_2(t)(I-A)+\cdots+f_k(t)(I-A)^{k-1}.
\end{array}$$ Obviously, this formula holds not only for all
$t\in[0,1]$, but also for any $t$ such that
$\sigma\left(I-(I-A)t\right)\cap\dR^-_0=\phi$. In particular, for
$t=1$ we may compute directly $\log A$: $$\log
A=f_1(1)I+f_2(1)(I-A)+\cdots+f_k(1)(I-A)^{k-1}.$$ }
\end{rem}

\section{Example}

To illustrate the method proposed, we consider the matrix
$$A=\left[\begin{array}{ccc} 7&4&-4\\4&7&-4\\-1&-1&4
\end{array}\right].$$ To compute $\log A$ we have to work
with the matrix $I-A$. The spectrum of $I-A$ is $\{-11,-2,-2\}$
and its minimum polynomial is $p(\lambda)=\lambda^2 +13\lambda
+22$. Applying directly (\ref{formula2}), we have
$$f_1(t)=\int_0^t\frac{22s}{1+13s+22s^2}\ ds,\quad
f_2(t)=\int_0^t\frac{-1}{1+13s+22s^2}\ ds.$$ Evaluating the
integrals, we may write
$$\log\left(I-(I-A)t\right)=f_1(t)I+f_2(t)(I-A),$$ where $t\in
[0,1]$ and $$\begin{array}{rcl}
f_1(t)&=&\frac{11}{9}\ln(1+2t)-\frac{2}{9}\ln(1+11t)\\f_2(t)&=&\frac{1}{9}
\ln\left(\frac{1+2t}{1+11t}\right).
\end{array}$$
Therefore
\begin{eqnarray*}
\log A&=&f_1(1)I+f_2(1)(I-A)\\&=&\left(\ln
3+\frac{2}{9}\ln\left(\frac{1}{4}\right)\right)I+\frac{1}{9}
\ln\left(\frac{1}{4}\right)(I-A).\end{eqnarray*}

{\bf Acknowledgements:} The author would like to thank Professor
F. Silva Leite for suggesting this problem and for useful
conversations on this topic.

\end{document}